\documentclass[12pt]{article} 
\usepackage[english]{babel}
\usepackage[latin1]{inputenc}
\usepackage{amsmath}
\usepackage{amsthm}
\usepackage{amssymb}
\usepackage{mathrsfs}
\usepackage{tikz-cd}
\usepackage{hyperref}
\usepackage{comment}
\usepackage{caption}

\usepackage{blkarray}
\usepackage{thmtools}

\usepackage[margin=0.6in]{geometry}

\newtheorem{theorem}{Theorem}[section]
\newtheorem{proposition}[theorem]{Proposition}
\newtheorem{corollary}[theorem]{Corollary}
\newtheorem{lemma}[theorem]{Lemma}

\theoremstyle{definition}

\newtheorem{definition}[theorem]{Definition}

\declaretheoremstyle[
numberlike=theorem,
bodyfont=\normalfont,
spaceabove=1em plus 0.75em minus 0.25em,
spacebelow=1em plus 0.75em minus 0.25em,
qed={$\diamondsuit$},
]{myExampleStyle}

\declaretheorem[
style=myExampleStyle,
title=Example,
refname={example,examples},
Refname={Example,Examples}
]{example}

\newcommand{\trop}{\operatorname{trop}}
\newcommand{\mult}{\operatorname{mult}}
\newcommand{\Pic}{\operatorname{Pic}}
\newcommand{\Hom}{\operatorname{Hom}}

\title{Tropical computations for toric intersection theory in Macaulay2}
\author{Alessio Borzì}
\date{}

\begin{document}
	
	\maketitle
	
	\begin{abstract}
		We present the Macaulay2 package \verb|TropicalToric.m2| for toric intersection theory computations using tropical geometry.
	\end{abstract}
	
	\section{Introduction}
	
	
	Toric varieties are ubiquitous in algebraic geometry. Their intersection theory was first studied in Fulton and Sturmfels \cite{fulton1997toricintersection}, and it has many applications in different contexts, including: wonderful and tropical compactifications \cite{tevelev2007compactification, deconcini1995wonderful, feichtner2004chow}, birational geometry \cite{gibeny2010equations, gibeny2012bounds, castravet2020blown}, tropical geometry \cite{katz2011realization}, tropical intersection theory \cite{katz2009toolkit, osserman2013lifting, allerman2010tropicalintersection, katz2012tropicalfromtoric, shaw2013matroidal} and combinatorial Hodge theory \cite{huh2012milnor,huhkatz,adiprasito2018hodge}.
	
	In a certain way, intersection classes of a toric variety with fan $\Sigma$ can be thought in terms of balanced subfans of $\Sigma$, also referred as \emph{Minkowski weights} (see \cite{fulton1997toricintersection} or \cite[Theorem 6.7.5]{maclagan2015tropical}). From the Structure Theorem of Tropical Geometry \cite[Theorem 3.3.5]{maclagan2015tropical} we know that the tropicalization of a subvariety of a torus $Y \subseteq T^n$ is a balanced fan. A surprising connection between tropical and toric geometry is that the tropicalization of $Y$ is the balanced fan corresponding to the intersection class of the closure of $Y$ inside an ``enough refined" toric variety (see Theorem \ref{thm:the theorem of the algorithm} for a more precise statement). This fact allows us to compute toric intersection classes starting from the data of the tropicalization.
	
	We present a new package \verb|TropicalToric.m2| for Macaulay2 \cite{M2} available at the following link.
	\begin{center}
		\url{https://github.com/AlessioBorzi/TropicalToric}
	\end{center}
	The package implements toric cycles and intersection products on simplicial toric varieties (Section \ref{sec2}), and, following the ideas outlined above, allows us to compute the intersection class of an irreducible subvariety of a simplicial toric variety not contained in the toric boundary, from the data of its tropicalization (Section \ref{sec3}). The tropicalization is performed with the use of the Macaulay2 package \verb|Tropical.m2| \cite{amendola2017computing}. Further, we present some applications to the intersection theory of wonderful compactifications and the moduli space $\overline{M}_{0,n}$, and illustrate an example in a multiprojective space using a theorem of Huh and Katz \cite{huhkatz} about characteristic polynomials of realizable matroids.
	
	\paragraph{Acknowledgements} I am grateful to my advisor, Diane Maclagan, for her guidance and suggestions during the development of this project. Part of the Macaulay2 package presented is based on some code previously written by Diane Maclagan, Sameera
	Vemulapalli, Corey Harris, Erika Pirnes and Ritvik Ramkumar.
	
	\section{Toric intersection theory}\label{sec2}
	
	
	In this section, we review the basics of toric intersection theory, for more information see \cite{fulton1997toricintersection}, \cite[Section 12.5]{cox2011toric} or \cite[Section 6.7]{maclagan2015tropical}. In addition, we showcase how it is implemented in the package.
	
	Let $X_\Sigma$ be a smooth complete toric variety of dimension $n$ with fan $\Sigma$. We denote by $\Sigma(k)$ the cones of $\Sigma$ of dimension $k$, with $Z^k(X_\Sigma) = Z_{n-k}(X_\Sigma)$ the group of codimension-$k$ cycles and with $A^k(X_\Sigma) = A_{n-k}(X_\Sigma)$ the codimension-$k$ Chow group, that is, the group of codimension-$k$ cycles modulo rational equivalence. The codimension-$k$ Chow group $A^k(X_\Sigma)$ is generated by the set $\{ [V(\sigma)] : \sigma \in \Sigma(k) \}$ of classes of orbit closures of codimension-$k$, and the relations in each Chow group can be described in an explicit way, see \cite[Proposition 2.1]{fulton1997toricintersection}.
	
	We can show (see \cite[Chapter 8]{fulton1998intersection}) that there is an intersection product $A^k(X_\Sigma) \times A^r(X_\Sigma) \rightarrow A^{k+r}(X_\Sigma)$ that makes $A^*(X_\Sigma) = \bigoplus_{k=0}^n A^k(X_\Sigma)$ into a graded ring, called the \emph{Chow ring} of $X_\Sigma$. If we now assume that $X_\Sigma$ is just complete and simplicial, then the intersection product can be defined on rational cycles, making $A^*(X_\Sigma)_{\mathbb{Q}} = A^*(X_\Sigma) \otimes \mathbb{Q}$ into a graded ring. The structure of the Chow ring has an explicit description, see for instance \cite[Theorem 6.7.1]{maclagan2015tropical}.

	Our Macaulay2 package implements toric cycles and the intersection product as in \cite[Lemma 12.5.2]{cox2011toric}.
	
	\begin{example}
		Let $X_\Sigma$ be the blow-up of $\mathbb{P}^2$ at one of the coordinate points, where the fan $\Sigma$ and the first lattice points of its rays are as in Figure \ref{fig1}.
		\begin{figure}[h!]
			\centering
			\definecolor{azure}{RGB}{80,150,255}
			\begin{tikzpicture}[thick,scale=1,baseline=(current bounding box.center)]
			\coordinate (O) at (0,0);
			
			\coordinate (V1) at (1,0);
			\coordinate (V2) at (0,1);
			\coordinate (V3) at (-1,0);
			\coordinate (V4) at (0,-1);
			
			\coordinate (A1) at (1,1);
			\coordinate (A2) at (-1,1);
			\coordinate (A3) at (-1,-1);
			\coordinate (A4) at (1,-1);
			
			\begin{scope}[thick,dashed,,opacity=0.6]
			\draw (O) -- (V1);
			\draw (O) -- (V2);
			\draw (O) -- (V3);
			\draw (O) -- (V4);
			\end{scope}
			
			\fill[color=orange,opacity=0.6] (O) -- (V1) -- (A1);
			\fill[color=azure,opacity=0.6] (O) -- (V2) -- (A1);
			\fill[color=yellow, opacity=0.6] (V2) -- (A2) -- (A3) -- (O);
			\fill[color=green, opacity=0.6] (V1) -- (A4) -- (A3) -- (O);
			
			\draw[line width=1] (O) -- (V1);
			\draw[line width=1] (O) -- (V2);
			\draw[line width=1] (O) -- (A3);
			\draw[line width=1] (O) -- (A1);
			
			\node at (1.3,0) {$\rho_0$};
			\node at (1.3,1.3) {$\rho_1$};
			\node at (0,1.3) {$\rho_2$};
			\node at (-1.3,-1.3) {$\rho_3$};
			
			
			\filldraw (0,0) circle (0.025);
			\end{tikzpicture}
			\hspace{2cm}
			\parbox[t][0cm][c]{4cm}{$\begin{blockarray}{cccc}
				\rho_0 & \rho_1 & \rho_2 & \rho_3 \\
				\begin{block}{[cccc]}
				1 & 1 & 0 & -1 \\
				0 & 1 & 1 & -1 \\
				\end{block}
			\end{blockarray}$}
		\captionsetup{width=0.8\textwidth}
		\caption{}\label{fig1}
		\end{figure}
		
		\noindent
		Let $H$ be the strict transform in $X_\Sigma$ of a general line in $\mathbb{P}^2$ and $E$ be the exceptional divisor. The Picard group of $X_\Sigma$ is generated by the classes of these two divisors $\Pic(X_\Sigma) = \langle [H],[E] \rangle$. With the notation above, we have
		\begin{align*}
		[V(\rho_0)] &= [H] - [E] & [V(\rho_1)] &= [E] \\
		[V(\rho_2)] &= [H] - [E] & [V(\rho_3)] &= [H]
		\end{align*}
		Now we verify with our package that the divisor class $[V(\rho_1)]$ has a negative self-intersection.
		\begin{verbatim}
	i1 : needsPackage "TropicalToric";
	i2 : raysList = {{1,0},{1,1},{0,1},{-1,-1}};
	i3 : coneList = {{0,1},{1,2},{2,3},{3,0}};
	i4 : X = normalToricVariety (raysList, coneList);
		\end{verbatim}
		Now define the toric cycle $V(\rho_1)$.
		\begin{verbatim}
	i5 : E = X_{1}
	o5 = X
	      {1}
	o5 : ToricCycle on X
		\end{verbatim}
		The type \verb|ToricCycle| should not be confused with the type \verb|ToricDivisor| from the \verb|NormalToricVarieties| package. The toric cycle $V(\sigma)$ of the normal toric variety \verb|X| associated to the cone $\sigma$ given by a list of rays \verb|L| is defined with the command \verb|X_L|. For example, \verb|X_{1,2}| or \verb|X_{0}| define toric cycles, whereas \verb|X_1| defines a toric divisor. We are allowed only to multiply a toric cycle with a toric divisor. Now, we finally compute the self intersection of $E$:
		\begin{verbatim}
	i6 : X_1 * E
	o6 = - X
	        {1, 2}
	o6 : ToricCycle on X
		\end{verbatim}
		The resulting cycle $-V(\rho_1+\rho_2)$ is rationally equivalent to $E^2$. The negative sign tells us that the self-intersection number of the exceptional divisor is $-1$. We can compute the degree of maximal codimension cycles with \verb|degCycle|:
		\begin{verbatim}
	i7 : degCycle(-X_{1,2})
	o7 = -1
		\end{verbatim}
	\vspace{-0.8cm}\qedhere

	\end{example}
	
	\section{Tropical computations}\label{sec3}
	
	In this section we describe and showcase the algorithm implemented in the main function of the package \verb|classFromTropical| that computes the intersection class of an irreducible subvariety of a simplicial toric variety from its tropicalization.
	
	The algorithm is mainly based on the following result, that appears in various versions in the literature, see for instance {\cite[Lemma 2.3]{katz2011realization}, \cite[Section 9]{katz2009toolkit} or \cite[Theorem 6.7.7]{maclagan2015tropical}}.
	
	\begin{theorem}\label{thm:the theorem of the algorithm}
		Let $Y$ be a subvariety of the algebraic torus $T^n$ and let $\overline{Y}$ be its closure in a toric variety $X_\Sigma$ such that $|\Sigma| = \trop(Y)$ and $\Sigma$ is simplicial. Let $\Sigma'$ be a completion of the fan $\Sigma$ and let $i:X_{\Sigma} \rightarrow X_{\Sigma'}$ be the induced inclusion. Then, for every maximal cone $\sigma$ in $\Sigma$ we have
		\[ m(\sigma) = \deg \big( [i_*(\overline{Y})] \cdot [V(\sigma)] \big), \]
		where $m(\sigma)$ is the multiplicity of $\sigma$ in $\trop(Y)$.
	\end{theorem}
	
	
	Now let $Y$ be an irreducible $k$-dimensional subvariety of an $n$-dimensional simplicial toric variety $X_\Sigma$, and suppose that $Y \cap T^n \neq \emptyset$. Note that in this setting we cannot directly apply Theorem \ref{thm:the theorem of the algorithm} since $\trop(Y)$ is not necessarily a subfan of $\Sigma$.
	
	In order to compute the class of $Y$ in the Chow ring of $X_\Sigma$, we proceed as follows. First, let $\Sigma'$ be a completion of $\Sigma$ and $i:X_\Sigma \rightarrow X_{\Sigma'}$ be the induced inclusion. Now let $\widetilde{\Sigma}$ be a refinement of $\Sigma'$ such that it contains a subfan with support the tropicalization of $Y \cap T^n$, and let $\pi: X_{\widetilde{\Sigma}}	\rightarrow X_{\Sigma'}$ be the induced toric map. From \cite[Proposition 2.4]{fulton1997toricintersection} we have an isomorphism $A_k(X_{\Sigma'}) \simeq \Hom(A^k(X_{\Sigma'}),\mathbb{Z})$ mapping a class $[Z]$ to the homomorphism $[Z'] \mapsto \deg([Z]\cdot [Z'])$. Therefore, in order to compute the class $[Y] \in A_k(X_{\Sigma})$, it is enough to compute the intersection numbers $\deg([i_*(Y)] \cdot [V(\sigma)])$ for every $\sigma \in \Sigma'(k)$, as the classes $[V(\sigma)]$ generate $A^k(X_{\Sigma'})$. Let $Y'$ be the strict transform of $Y$ in $X_{\widetilde{\Sigma}}$. From the projection formula \cite[Proposition 2.3 (c)]{fulton1998intersection}, we have
	\[ [i_*(Y)] \cdot [V(\sigma)] = \pi_*([Y']) \cdot [V(\sigma)] = \pi_*\big( [Y'] \cdot \pi^*([V(\sigma)]) \big) \]
	from which it follows $\deg \big([i_*(Y)] \cdot [V(\sigma)] \big) = \deg \big([Y'] \cdot \pi^*([V(\sigma)]) \big)$. These last intersection numbers can be computed from the tropicalization of $Y \cap T^n$ by using Theorem \ref{thm:the theorem of the algorithm}, since $\deg \big( [Y'] \cdot [V(\sigma')] \big)$ is the multiplicity of the cone $\sigma' \in \widetilde{\Sigma}(k)$ in the tropicalization of $Y \cap T^n$.
	
	The function \verb|classFromTropical| performs the above algorithm to compute a toric cycle rationally equivalent to a given irreducible subvariety $Y$ of a simplicial toric variety $X_\Sigma$. The input of the function consists of the toric variety $X_\Sigma$ and the ideal $I$ of $Y \cap T^n$ of the Laurent ring of $T^n$. Since Laurent rings are not implemented in Macaulay2, the actual input will be instead the saturation of $I$ with respect to the product of the variables in the polynomial ring:
	\begin{verbatim}
i2 : X = toricProjectiveSpace 2;
i3 : R = QQ[x,y];
i4 : I = ideal(x+y+1);
i5 : classFromTropical(X,I)
o5 = X
      {0}
o5 : ToricCycle on X
i6 : J = ideal(x*y + x + y);
i7 : classFromTropical(X,J)
o7 = 2*X
        {0}
o7 : ToricCycle on X
	\end{verbatim}
	\noindent
	The function \verb|classFromTropicalCox| allows us to input the ideal of $Y$ in the Cox ring of $X_\Sigma$: 
	
	\begin{verbatim}
i8 : R = ring X;
i9 : I = ideal(R_0+R_1+R_2);
i10 : classFromTropicalCox(X,I)
o10 = X
       {0}
o10 : ToricCycle on X
	\end{verbatim}
	
	\section{Applications}\label{sec4}

	\subsection{Wonderful compactifications}\label{sec:wonderful}

	
	Let $\mathcal{A}$ be an essential hyperplane arrangement of $n+1$ hyperplanes in $\mathbb{P}^d$. The intersection lattice $\mathcal{L}(\mathcal{A})$ of $\mathcal{A}$ is isomorphic to the lattice of flats of the underlying matroid $M$ of $\mathcal{A}$ \cite[Proposition 3.6]{stanley2004introduction}. Fix a building set $\mathcal{G}$ of the lattice of flats of $M$ (see \cite[Section 2]{feichtner2004chow}), let $\Sigma \subseteq \mathbb{R}^{n+1} / \mathbb{R} \textbf{1} \simeq \mathbb{R}^n$ be the Bergman fan of $M$ with respect to $\mathcal{G}$ (see \cite[Chapter 4]{maclagan2015tropical}) and let $X_{\Sigma}$ be its associated toric variety. From \cite[Proposition 4.1.1]{maclagan2015tropical} the hyperplane arrangement complement $Y = \mathbb{P}^d \setminus \cup \mathcal{A}$ is naturally isomorphic to a linear subspace of the algebraic torus $T^n$. Thus we can embed $Y$ inside the toric variety $X_{\Sigma}$ and consider its closure $\overline{Y}$. This compactification coincides with the so-called De Concini-Procesi \emph{wonderful compactification} \cite{deconcini1995wonderful}, with respect to the building set $\mathcal{G}$ (see \cite[Section 4]{tevelev2007compactification}). The next result follows from \cite[Theorem 3.2]{deconcini1995wonderful} (see also \cite[Definition 2.3]{feichtner2005deconciniprocesi}).
	
	\begin{proposition}\label{prop:compactification from blow-up}
		Let $X_1,\dots,X_t$ be a linear extension of the opposite order of $\mathcal{L}(\mathcal{A})$. The wonderful compactification $\overline{Y}$ is the result of successively blowing up $\mathbb{P}^d$ at (the strict transforms of) $X_1,\dots,X_t$.
	\end{proposition}
	
	In \cite{feichtner2004chow} Feichtner and Yuzvinsky showed that the cohomology of $\overline{Y}$ agrees with that of $X_{\Sigma}$. Since both varieties are smooth, the Chow ring is the same object as the cohomology ring, and we obtain the following result.
	
	\begin{theorem}[{\cite[Theorem 6.7.14]{maclagan2015tropical}}]\label{thm:chow ring}
		Let $\overline{Y}$ be a wonderful compactification of a hyperplane arrangement $\mathcal{A}$ with respect to a building set $\mathcal{G}$, and let $\Sigma$ be the associated Bergman fan. Then
		\[ A^*(\overline{Y}) \simeq A^*(X_{\Sigma}). \]
	\end{theorem}
	
	The above theorem allows us to view the intersection classes of a wonderful compactification as intersection classes of the associated toric variety. Thus, we can use our package to perform intersection theory computations on wonderful compactifications.
	
	\begin{example}\label{ex:line arrangement}
		Let $\mathcal{A}$ be a line arrangement consisting of $4$ lines $L_0,L_1,L_2,L_3$ in $\mathbb{P}^2$ given by the equations $x_0 = 0$, $x_1=0$, $x_2=0$, $x_0+x_1=0$ respectively. Let $A$ be the matrix with columns the normal vectors of the lines $L_i$, and let $P_1,P_2,P_3,P_4$ be the points of intersection of the lines of $\mathcal{A}$ as in the figure below.
		\begin{figure}[h!]
			\centering
			
			\parbox[t][-4cm][c]{5cm}{$A = \begin{bmatrix}
				1 & 0 & 0 & 1 \\
				0 & 1 & 0 & 1 \\
				0 & 0 & 1 & 0 \\
				\end{bmatrix}$}
			\begin{tikzpicture}
			[scale=1,auto=center,every node/.style={}]
			
			\draw[fill=black] (0,0) circle (0.08);
			\draw[fill=black] (1,0) circle (0.08);
			\draw[fill=black] (2,0) circle (0.08);
			\draw[fill=black] (1,1) circle (0.08);
			
			\coordinate (A1) at (0,2);
			\coordinate (A2) at (1,2.5);
			\coordinate (A3) at (2,2);
			\coordinate (A4) at (-1,-1);
			\coordinate (A5) at (1,-1.5);
			\coordinate (A6) at (3,-1);
			\coordinate (A7) at (-1,0);
			\coordinate (A8) at (3,0);
			
			\draw (-0.3,1.7) node {$L_0$};
			\draw (0.65,-1) node {$L_3$};
			\draw (2.3,1.7) node {$L_1$};
			\draw (3.4,0) node {$L_2$};
			
			\draw (-0.3,0.3) node {$P_3$};
			\draw (1.3,-0.3) node {$P_4$};
			\draw (2.3,0.3) node {$P_2$};
			\draw (1.5,1) node {$P_1$};
			
			\foreach \from/\to in
			{A1/A6,A2/A5,A3/A4,A7/A8}
			\draw (\from) -- (\to);
			
			\end{tikzpicture}
		\end{figure}
		
		\noindent
		The underlying matroid $M$ of $\mathcal{A}$, on the ground set $\{ 0,1,2,3 \}$, is realized by the matrix $A$ by labeling the columns with $0,1,2,3$ respectively. The lattice of flats $\mathcal{L}(M)$ of $M$ is represented by the diagram in Figure \ref{fig:2}.
		\begin{figure}[h!]
			\centering
			\begin{tikzpicture}[scale=1, align=center]
			\node (00) at (0,0) {$\emptyset$};
			
			\node (01) at (-1.5,1) {$\{0\}$};
			\node (11) at (-0.5,1) {$\{1\}$};
			\node (21) at (0.5,1) {$\{2\}$};
			\node (31) at (1.5,1) {$\{3\}$};
			
			\node (02) at (-3,2.5) {$\{0,1,3\}$};
			\node (12) at (-1,2.5) {$\{0,2\}$};
			\node (22) at (1,2.5) {$\{1,2\}$};
			\node (32) at (3,2.5) {$\{2,3\}$};
			
			\node (03) at (0,3.5) {$\{0,1,2,3\}$};

			\draw (00) -- (01);
			\draw (00) -- (11);
			\draw (00) -- (21);
			\draw (00) -- (31);
			
			\draw (01) -- (02);
			\draw (11) -- (02);
			\draw (31) -- (02);
			
			\draw (01) -- (12);
			\draw (21) -- (12);
			
			\draw (11) -- (22);
			\draw (21) -- (22);
			
			\draw (21) -- (32);
			\draw (31) -- (32);
			
			\draw (03) -- (02);
			\draw (03) -- (12);
			\draw (03) -- (22);
			\draw (03) -- (32);
			\end{tikzpicture}
			\caption{}\label{fig:2}
		\end{figure}
		
		\noindent
		There are $4$ rank $1$ flats, corresponding to the lines $L_0,L_1,L_2,L_3$, and $4$ rank $2$ flats, corresponding to the points $P_1,P_2,P_3,P_4$. Let $\mathcal{G} = \mathcal{L}(M) \setminus \{ \emptyset \}$ be the maximal building set of $\mathcal{L}(M)$. Then, from Proposition \ref{prop:compactification from blow-up}, the wonderful compactification $\overline{Y}$ of the complement $Y = \mathbb{P}^2 \setminus \cup \mathcal{A}$ with respect to $\mathcal{G}$ is the blow-up of $\mathbb{P}^2$ at the points $P_1,P_2,P_3,P_4$.
		In particular $\overline{Y}$ is a smooth projective surface, all Weil divisors are Cartier \cite[Proposition II.6.11]{hartshorne1977algebraicgeometry}, and the class group is isomorphic to the Picard group \cite[Corollary II.6.16]{hartshorne1977algebraicgeometry}. From \cite[Proposition V.3.2]{hartshorne1977algebraicgeometry}, the Picard group of $\overline{Y}$ has a basis given by
		\begin{equation}\label{eq: Picard basis}
			\Pic(\overline{Y}) = \langle [H], [E_1], \dots, [E_t] \rangle,
		\end{equation}
		where $[H]$ is the class of the strict transform $H$ of a general line in $\mathbb{P}^2$, and $[E_i]$ is the class of the exceptional divisor $E_i$ of the blow-up at $P_i$.

		The Bergman fan $\Sigma \subseteq \mathbb{R}^4 / \mathbb{R} \textbf{1}$ of $M$ with respect to $\mathcal{G}$ has $8$ rays, denoted $\{ \rho_i : 0 \leq i \leq 7 \}$. Their first lattice points are given by the columns of the following matrix
		\[ \begin{blockarray}{cccccccc}
		\rho _0 & \rho_1 & \rho_2 & \rho_3 & \rho_4 & \rho_5 & \rho_6 & \rho_7 \\
		\begin{block}{[cccccccc]}
		1 & 0 & 0 & 0 & 1 & 1 & 0 & 0 \\
		0 & 1 & 0 & 0 & 1 & 0 & 1 & 0 \\
		0 & 0 & 1 & 0 & 0 & 1 & 1 & 1 \\
		0 & 0 & 0 & 1 & 1 & 0 & 0 & 1 \\
		\end{block} 
		\end{blockarray} \]
		where $\rho_0,\rho_1,\rho_2,\rho_3$ correspond to the rank $1$ flats in $\mathcal{G}$, which in turn correspond to the lines $L_0,L_1,L_2,L_3$ respectively, and $\rho_4,\rho_5,\rho_6,\rho_7$ correspond to the rank $2$ flats in $\mathcal{G}$, that correspond to the points $P_1,P_2,P_3,P_4$ respectively. Since $\mathcal{G}$ is the maximal building set, the maximal cones of $\Sigma$ are just the maximal chains of the lattice of flats $\mathcal{L}(M)$.
		
		By using the isomorphism in Theorem \ref{thm:chow ring}, let $[Y_{\rho_i}]$ denote the class in $A^*(\overline{Y})$ isomorphic to the class of the torus invariant divisor of $X_\Sigma$ associated to the ray $\rho_i$. Expressing these divisors in the Picard basis \eqref{eq: Picard basis} we have:
		\begin{align}\label{eq:example rays}
		[Y_{\rho_0}] &= [H] - [E_1] - [E_2] & [Y_{\rho_4}] &= [E_1] \\
		[Y_{\rho_1}] &= [H] - [E_1] - [E_3] & [Y_{\rho_5}] &= [E_2] \nonumber \\
		[Y_{\rho_2}] &= [H] - [E_2] - [E_3] - [E_4] & [Y_{\rho_6}] &= [E_3] \nonumber \\
		[Y_{\rho_3}] &= [H] - [E_1] - [E_4] & [Y_{\rho_7}] &= [E_4] \nonumber
		\end{align}
		Let $\mathbb{C}[y_0^{\pm 1},y_1^{\pm 1},y_2^{\pm 1}]$ be the Laurent ring of the torus
		\[ T^3 = \{ (1:y_0:y_1:y_2) : y_0,y_1,y_2 \in \mathbb{C}^* \} \subseteq \mathbb{P}^3. \]
		The embedding $Y \hookrightarrow T^3$ is given by $(x_0:x_1:x_2) \mapsto (x_0:x_1:x_2:x_0+x_1)$, and the Laurent ideal of $Y$ inside $T^3$ is $I = (-1-y_0+y_2)$.
		
		Now let $C$ be the conic in $\mathbb{P}^2$ passing through $P_1$, $P_2$ and $P_3$ given by the equation $x_0x_1+x_0x_2+x_1x_2$. The ideal of $C$ in $T^3$ is $(y_0+y_1+y_0y_1) + I$. We expect the class of its strict transform in $\overline{Y}$ to be $[2H - E_1 - E_2 - E_3]$. We now verify this with our package, using the function \verb|classWonderfulCompactification|:
		\begin{verbatim}
i2 : R = QQ[y_0,y_1,y_2];
i3 : I = ideal(-1-y_0+y_2);
i4 : f = y_0+y_1+y_0*y_1;
i5 : raysList = {{-1,-1,-1},{1,0,0},{0,1,0},
                 {0,0,1},{0,-1,0},{-1,0,-1},
                 {1,1,0},{0,1,1}};
i6 : conesList = {{4,0},{4,1},{4,3},{5,0},{5,2},
                  {6,1},{6,2},{7,2},{7,3}};
i7 : X = normalToricVariety (raysList, conesList);
i8 : D = classWonderfulCompactification(X,I,f)
o8 = X    + X    + X
      {0}    {4}    {1}
o8 : ToricCycle on X
		\end{verbatim}
	To check that this is the result we expect, compare with \eqref{eq:example rays}. Note that we have (tropically) dehomogenized the rays of $X_\Sigma$ with respect to the first coordinate in order to be consistent with our choice of coordinates of $T^3$.
	\end{example}
	
	
	\subsection{The moduli space $\overline{M}_{0,n}$}
	
	The Deligne-Mumford compactification of the moduli space $M_{0,n}$ can be realized as a wonderful compactification (see for instance \cite[Example 6.7.16]{maclagan2015tropical}). Therefore, we can apply to $\overline{M}_{0,n}$ the machinery described in the previous section. As an application, we compute one of the $15$ Keel-Vermeire divisors of $\overline{M}_{0,6}$, using one of the equations listed in \cite[Table 2]{guillen2017presentation}. These divisors, found independently by Keel and Vermeire \cite{vermeire2002counterexample}, were the first example of an effective divisor of $\overline{M}_{0,n}$ which class lies outside the cone generated by the classes of the boundary divisors, answering in the negative to a conjecture of Fulton (see \cite{keel2013contractible}). 
	
	\begin{verbatim}
i2 : R = QQ[x_0..x_8];
i3 : I = ideal {-x_0+x_3+x_4, -x_1+x_3+x_5,-x_2+x_3+x_6,
                -x_0+x_2+x_7, -x_1+x_2+x_8, -x_0+x_1+1};
i4 : X = normalToricVariety fan tropicalVariety I;
i5 : f = x_0*x_1-x_2*x_3;
i6 : D = classWonderfulCompactification(X,I,f);
i7 : D = toricDivisorFromCycle(D)
o7 = X  - X  - 2*X  + X  + 2*X  + 2*X   - X   + 2*X   + 2*X   - X
      2    5      6    7      9      10    11      13      14    17
o7 : ToricDivisor on X
	\end{verbatim}
	Now fix the Picard basis of $X_\Sigma$ given by the boundary divisors associated to the rays of $\Sigma$ with first lattice points not equal to the standard vectors $e_i$. The complement of this Picard basis is indexed by the list $l = \{ 0,1,2,4,5,7,11,13,21 \}$. The function \verb|makeTransverse| computes a divisor linearly equivalent to a given divisor \verb|D|, with support disjoint from a given list \verb|l|. We use this function to compute a representation of the class of the Keel-Vermeire divisor computed above, in the Picard basis we fixed:
	\begin{verbatim}
i8 : l = {0,1,2,4,5,7,11,13,21};
i9 : D = makeTransverse(D,l)
o9 = X  - X  - X  + X  + 2*X   - X   - X   + X   + 2*X   + 2*X   
      3    6    8    9      10    14    16    17      18      19
      + 2*X   - X   - X
           20    22    24
o9 : ToricDivisor on X
	\end{verbatim}
	Finally, we verify that the obtained divisor is outside the cone generated by the classes of boundary divisors. In order to do so, we interface with Polymake \cite{polymake:2000} by using the function \verb|polymakeConeContains|:
	\begin{verbatim}
i10 : D = apply(#rays X, i->D#i);
i11 : Bdivisors = apply(#rays X, i-> makeTransverse(X_i,l));
i12 : Bdivisors = apply(Bdivisors, B-> apply(#rays X, i->B#i));
i13 : polymakeConeContains(D,Bdivisors)
o13 = false
	\end{verbatim}
	
	In \cite{hassett2002effective} it was proved, by using computational methods, that the boundary divisors and the Keel-Vermeire divisors generate the effective cone of $\overline{M}_{0,6}$. In \cite{castravet2020blown} it was proved that the effective cone of $\overline{M}_{0,n}$ for $n \geq 10$ is not polyhedral. The problem of determining the effective cone of $\overline{M}_{0,n}$ for $n \in \{ 7,8,9 \}$ is still open. Some examples of extremal effective divisors on $\overline{M}_{0,7}$ were found in \cite{castravet2013hypertree,opie2016extremal,doran2017simplicial}. We performed computations similar to those displayed above on $\overline{M}_{0,7}$ and found the mentioned examples with a brute-force approach. More recently in \cite{sikiric2022extreme} several thousands of extremal effective divisors on $\overline{M}_{0,7}$ were found.
	
	
	\subsection{Characteristic polynomials}
	
	Our last application is an explicit verification of a theorem proved by Huh and Katz \cite{huhkatz} about characteristic polynomials of realizable matroids.
	
	Let $\mathcal{A}$ be an arrangement of $n+1$ hyperplanes on $\mathbb{P}^d$, let $M$ be its underlying matroid of rank $d+1$, and let $\mathcal{L}(M)$ be the lattice of flats of $M$. The \emph{characteristic polynomial} of $M$ is
	\[ \chi_M(q) = \sum_{F \in \mathcal{L}(M)} \mu(\emptyset,F) q^{d+1-r(F)} \]
	where $\mu$ is the M\"{o}bius function of $\mathcal{L}(M)$ (see \cite[Section 3.7]{stanley2021enumerativeI}). The \emph{reduced characteristic polynomial} of $M$ is $\overline{\chi}_M(q) = \chi_M(q)/(q-1)$.
	
	Now we embed the complement $Y = \mathbb{P}^{d} \setminus \cup \mathcal{A}$ in $T^n \subseteq \mathbb{P}^n$, as described in Section \ref{sec:wonderful}, and consider the Cremona map
	\[ \varphi: \mathbb{P}^n \dashrightarrow \mathbb{P}^n, \quad (x_0,\dots,x_n) \mapsto (x_0^{-1},\dots,x_n^{-1}).  \]
	Finally, let $\overline{Z}$ be the closure in $\mathbb{P}^n \times \mathbb{P}^n$ of the graph $Z$ of the restriction $\varphi_{|Y}$.
	
	\begin{theorem}[Huh-Katz \cite{huhkatz}]
		Define the integers $a_i \in \mathbb{Z}$ by the formula
		\[ \overline{\chi}_M(q) = \sum_{i=0}^d (-1)^i a_i q^{d-i}. \]
		Then
		\[ [\overline{Z}] = \sum_{i=0}^d a_i [\mathbb{P}^{r-i} \times \mathbb{P}^i] \in A_{d}(\mathbb{P}^n \times \mathbb{P}^n). \]
	\end{theorem}
	
	\begin{example}
		Let $G$ be the graph as in the figure below.
		
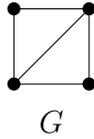
\begin{figure}[h!]
			\centering
			\begin{tikzpicture}
				\coordinate (A1) at (0,0);
				\coordinate (A2) at (1,0);
				\coordinate (A3) at (1,1);
				\coordinate (A4) at (0,1);
				
				\draw[fill=black] (A1) circle (0.08);
				\draw[fill=black] (A2) circle (0.08);
				\draw[fill=black] (A3) circle (0.08);
				\draw[fill=black] (A4) circle (0.08);
				
				\draw (A1) -- (A2);
				\draw (A1) -- (A3);
				\draw (A1) -- (A4);
				\draw (A2) -- (A3);
				\draw (A3) -- (A4);
				
				\node at (0.5,-0.5) {$G$};
			\end{tikzpicture}
			\hspace{2cm}
			\parbox[t][-2cm][c]{5cm}{$A = \begin{bmatrix}
				1 & 0 & 0 & 1 & 1 \\
				0 & 1 & 0 & -1& 0 \\
				0 & 0 & 1 & 0 & -1 \\
				\end{bmatrix}$}
		\end{figure}
	
		\noindent
		Let $M$ be the rank $3$ graphic matroid of $G$, realized by the matrix $A$ above. The characteristic polynomial of $M$ coincide with the chromatic polynomial of $G$. Let $\mathbb{C}[x_0^{\pm 1},x_1^{\pm 1},x_2^{\pm 1},x_3^{\pm 1}]$ be the Laurent ring of the torus
		\[ T^4 = \{ (1:x_0:x_1:x_2:x_3) : x_0,x_1,x_2,x_3 \in \mathbb{C}^* \} \subseteq \mathbb{P}^4. \]
		Let $\mathcal{A}$ be the hyperplane arrangement realizing $M$. More explicitly, the normal vectors of its hyperplanes are the columns of the matrix $A$. The Laurent ideal of the hyperplane arrangement complement $Y = \mathbb{P}^2 \setminus \cup \mathcal{A}$ embedded in $T^4 \subseteq \mathbb{P}^4$ is given by $I = (-1+x_0+x_2,-1+x_1+x_3)$.
		
		Now consider a copy of $T^4$ with Laurent ring $\mathbb{C}[x_4^{\pm 1},x_5^{\pm 1},x_6^{\pm 1},x_7^{\pm 1}]$. Let $Z \subseteq T^4 \times T^4$ be the graph of the Cremona map $\varphi: \mathbb{P}^4 \dashrightarrow \mathbb{P}^4$ restricted to $Y$. The ideal of $Z$ in the Laurent ring $\mathbb{C}[x_0^{\pm 1},\dots,x_7^{\pm 1}]$ of $T^4 \times T^4$ is generated by $I$ and the polynomials $x_ix_{i+4} - 1$ for $i \in \{ 0,1,2,3 \}$.
		\begin{verbatim}
i2 : R = QQ[x_0..x_7];
i3 : I = ideal(-1+x_0+x_2,-1+x_1+x_3,
               x_0*x_4-1,x_1*x_5-1,x_2*x_6-1,x_3*x_7-1);
o3 : Ideal of R
i4 : P4 = toricProjectiveSpace 4;
i5 : X = NormalToricVarieties$cartesianProduct(P4,P4);
i6 : D = classFromTropical(X,I)
o6 = 4*X              + 4*X                   + X
        {0,1,2,3,5,6}      {0,1,2,5,6,7}         {0,1,5,6,7,8}
o6 : ToricCycle on X
		\end{verbatim}
		We obtained $[\overline{Z}] = [\mathbb{P}^2 \times \mathbb{P}^0] + 4 [\mathbb{P}^1 \times \mathbb{P}^1] + 4 [\mathbb{P}^0 \times \mathbb{P}^2]$. We now verify that the coefficients of this class are the same, up to sign, to those of the (reduced) chromatic polynomial of $G$:
	\begin{verbatim}
i7 : needsPackage "Graphs";
i8 : G = graph({{0,1},{1,2},{2,3},{3,0},{0,2}});
i9 : p = chromaticPolynomial G
       4     3     2
o9 = x  - 5x  + 8x  - 4x
o9 : ZZ[x]
i10 : x = (ring p)_0;
i11 : p/(x-1)
       3     2
o11 = x  - 4x  + 4x
o11 : frac(ZZ[x])
	\end{verbatim}
	\vspace{-0.8cm}
	\end{example}
	

	\bibliographystyle{abbrv}
	\bibliography{Reference.bib}
	
 	\noindent
	{\scshape Alessio Borz\`{i}} \quad \texttt{Alessio.Borzi@warwick.ac.uk}\\
 	{\scshape  Mathematics Institute, University of Warwick, Coventry CV4 7AL, United Kingdom.}
	
\end{document}